\documentclass{amsart}

\usepackage{amssymb}
\usepackage{amsmath}
\usepackage{latexsym}
\usepackage{mathrsfs}
\usepackage{comment}
\usepackage{graphicx}
\usepackage{caption}
\usepackage{color}
\usepackage{hyperref}

\numberwithin{equation}{section}
\newtheorem{thm}[equation]{Theorem}
\newtheorem{prop}[equation]{Proposition}

\newtheorem{rem}[equation]{Remark}
\newtheorem{lem}[equation]{Lemma}
\newtheorem{corol}[equation]{Corollary}

\providecommand{\abs}[1]{\left\lvert#1\right\rvert}

\begin{document}

\title[Steklov eigenvalues]
 {Neumann to Steklov eigenvalues: asymptotic and  monotonicity results}

%----------Author 1
\author{Pier Domenico Lamberti}

\address{%
Dipartimento di Matematica\\
Universit\`a degli Studi di Padova\\
Via Trieste, 63\\
%P.O. Box 133\\
35126 Padova\\
Italy}

\email{lamberti@math.unipd.it}

%\thanks{This work was completed with the support of our
%\TeX-pert.}
%----------Author 2
\author{Luigi Provenzano}
\address{Dipartimento di Matematica\\
Universit\`a degli Studi di Padova\\
Via Trieste, 63\\
%P.O. Box 133\\
35126 Padova\\
Italy}
\email{proz@math.unipd.it}
%----------classification, keywords, date

%\date{November 05, 2013}
%----------additions

%%% ----------------------------------------------------------------------

\begin{abstract}
We consider the Steklov eigenvalues of the Laplace operator as limiting Neumann eigenvalues in a problem of mass concentration at the boundary of a ball. We discuss the asymptotic behavior of the Neumann eigenvalues and find explicit formulas for their derivatives at the limiting problem. We deduce that the Neumann eigenvalues have a monotone behavior in the limit and that  
Steklov eigenvalues locally minimize the Neumann eigenvalues.
\end{abstract}

%%% ----------------------------------------------------------------------
\maketitle
%%% ----------------------------------------------------------------------
%\tableofcontents

%%%%%%%%%%%%%%%%%%%%%%%%%%%%%%%%%%%%%%%%%%%%%%%%%%%%%%%%%%%%%%%%%%%%%%%%%%%%%%%%%%%%%%%%%%%%%%%%%%%%
%%%%%%%%%% INTRODUCTION
%%%%%%%%%%%%%%%%%%%%%%%%%%%%%%%%%%%%%%%%%%%%%%%%%%%%%%%%%%%%%%%%%%%%%%%%%%%%%%%%%%%%%%%%%%%%%%%%%%%%
%%%%%%%%%%

\section{Introduction}\label{intro}

Let   $B$  be the unit ball in $\mathbb R^N$, $N\geq 2$, centered at zero.  We consider the Steklov eigenvalue  problem for the Laplace operator

\begin{equation}\label{Ste}
\left\{\begin{array}{ll}
\Delta u =0 ,\ \ & {\rm in}\  B,\\
\frac{\partial u}{\partial \nu }=\lambda\rho u ,\ \ & {\rm on}\  \partial B,
\end{array}\right.
\end{equation}
in the unknowns $\lambda$ (the eigenvalue) and $u$ (the eigenfunction), where $\rho =M/\sigma_N$, $M>0$ is a fixed constant, and $\sigma_N$ denotes the surface measure of $\partial B$. 

As is well-known the eigenvalues of problem (\ref{Ste}) are given explicitly by the sequence 
\begin{equation}\label{intro1}
\lambda_l=\frac{l}{\rho},\ \ \ \ l\in {\mathbb{N}},
\end{equation}
and  the eigenfunctions corresponding to $\lambda_l$ are the homogeneous harmonic polynomials of degree $l$. In particular, the multiplicity of $\lambda_l$ is 
$(2l+N-2) (l+N-3)! /(l!(N-2)!)$,
and only $\lambda_0$ is simple, the corresponding eigenfunctions being the constant functions. See \cite{folland} for an introduction to the theory of 
harmonic polynomials. 

A classical reference for problem (\ref{Ste}) is
\cite{stek}. For a recent survey paper, we refer to \cite{gipo}; see  also \cite{lambertisteklov}, \cite{lapro} for related problems.

It is well-known that for $N=2$, problem (\ref{Ste}) provides the vibration modes of a free elastic membrane the total  mass of which  is  $M$ and is concentrated at the boundary with density $\rho$; see e.g., \cite{bandle}. 
As is pointed out in \cite{lapro}, such a boundary concentration phenomenon can be explained in any dimension $N\geq 2$ as follows.

For any $0< \varepsilon <1$, we define a `mass density'  $\rho_{\varepsilon}$ in the whole of $B$ by setting  

\begin{equation}\label{densitaintro}
\rho_{\varepsilon}(x)=\left\{
\begin{array}{ll}
\varepsilon,& {\rm if\ }|x|\le 1-\varepsilon ,\\
\frac{M-\varepsilon \omega_N (1-\varepsilon )^N}{\omega_N  (1-(1-\varepsilon)^N  ) },& {\rm if\ }1-\varepsilon <|x|<1,
\end{array}
\right.
\end{equation}
where $\omega_N=\sigma_N/N$ is the measure of the unit ball.  Note that for any $x\in B $  we have $\rho_{\varepsilon}(x)\to 0$ as $\varepsilon \to 0$, and   $\int_{B}\rho_{\varepsilon}dx=M$ for all $\varepsilon>0$, which means that the 
`total mass' M is fixed and concentrates at the boundary of $B$ as $\varepsilon \to 0$. Then we consider the following eigenvalue problem 
for the Laplace operator with Neumann boundary conditions
\begin{equation}\label{Neu}
\left\{\begin{array}{ll}
-\Delta u =\lambda \rho_{\varepsilon} u,\ \ & {\rm in}\  B,\\
\frac{\partial u}{\partial \nu }=0 ,\ \ & {\rm on}\  \partial B. 
\end{array}\right.
\end{equation}
We recall  that for $N=2$ problem (\ref{Neu}) provides the vibration modes of a free elastic membrane with mass density $\rho_{\varepsilon }$ and total mass $M$ (see e.g., \cite{cohi}). The eigenvalues of (\ref{Neu}) have finite multiplicity and  form a sequence
$$
\lambda_0(\varepsilon)<\lambda_1(\varepsilon)\leq\lambda_2(\varepsilon)\leq\cdots,
$$
depending on $\varepsilon$,
with $\lambda_0(\varepsilon)=0$. 

It is not difficult to prove that for any $l\in {\mathbb{N}}$
\begin{equation}\label{intro1,5}
	\lambda_l(\varepsilon)\to \lambda_l,\ \ {\rm as}\ \varepsilon \to 0,
\end{equation}
see \cite{arr}, \cite{lapro}.  (See also \cite{buosoprovenzano} for a detailed analysis of the analogue problem for the biharmonic operator.) 
Thus the Steklov problem can be considered as a limiting Neumann problem where the mass is concentrated at the boundary of the domain.  

In this paper we study the asymptotic behavior of $\lambda_l(\varepsilon)$ as $\varepsilon\rightarrow 0$. Namely, we prove that such eigenvalues are continuously differentiable  with respect to $\varepsilon$ for $\varepsilon \geq 0$ small enough, and that the following formula holds
\begin{equation}\label{intro2}
\lambda'_l(0)= \frac{2l\lambda_l}{3}+\frac{2\lambda^2_l}{N(2l+N)}.
\end{equation}
In particular,  for $l\ne 0$,  $\lambda'_l(0)>0$  hence $\lambda_l(\varepsilon )$ is strictly increasing  and  the Steklov eigenvalues $\lambda_l$ minimize the Neumann eigenvalues $\lambda_l(\varepsilon)$ for $\varepsilon $ small enough. 

\begin{figure}[ht!]
 \centering
    \includegraphics[width=0.6\textwidth]{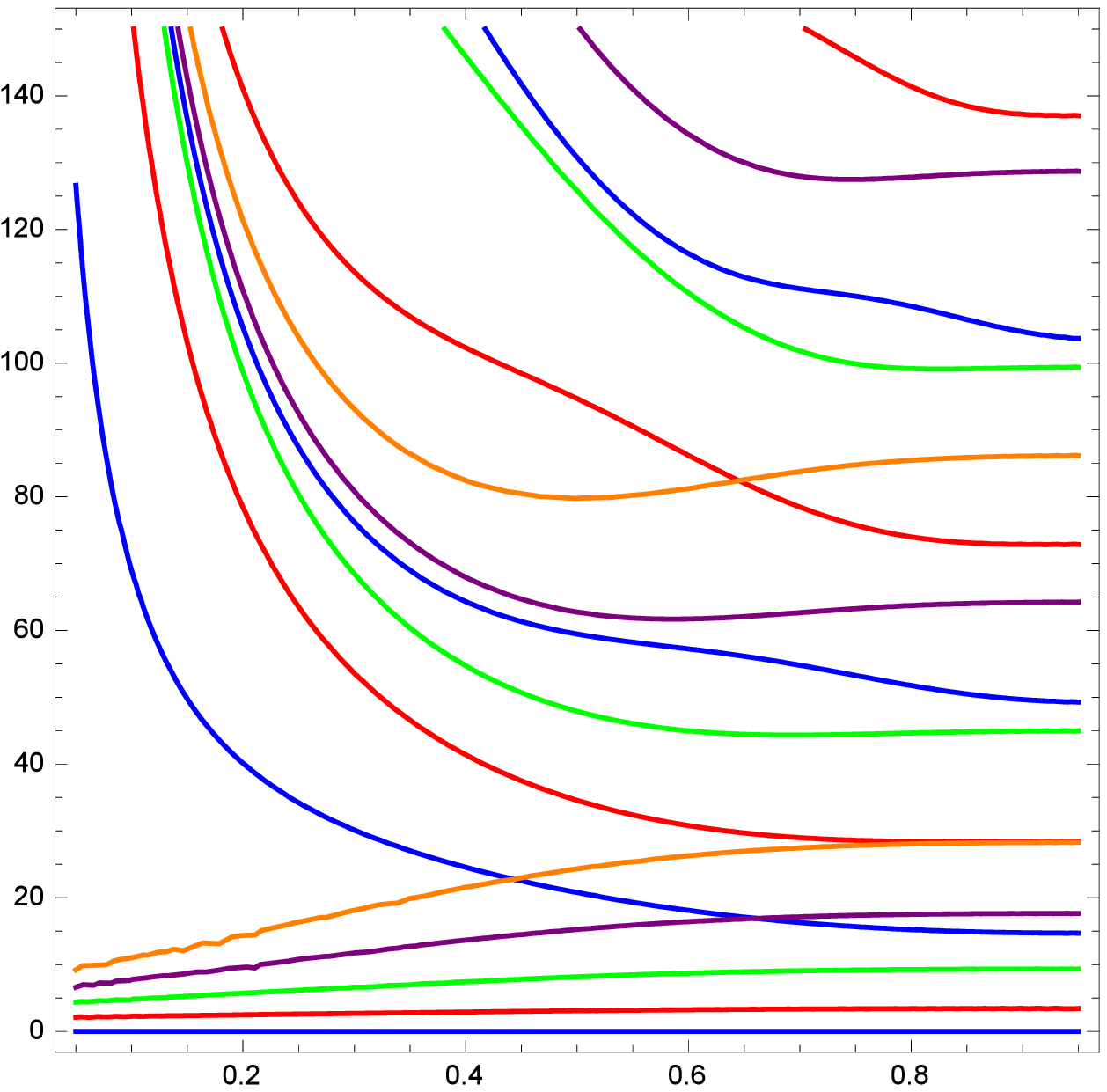}
		 \caption{{ Solution branches  of equation \eqref{implicitformula1N} with $N=2$, $M=\pi$ for $(\varepsilon,\lambda)\in]0,1[\times]0,150[$ . The colors refer to the choice of $l$ in \eqref{implicitformula1N}: blue ($l=0$), red ($l=1$), green ($l=2$), purple ($l=3$), orange ($l=4$).} 
	%		Branches of eigenvalues in the case $N=2$, $M=\pi$. The colors refer to the index $l\in\mathbb N$, in particular blue ($l=0$), red ($l=1$), green ($l=2$), purple ($l=3$), orange ($l=4$). The range considered is $(\varepsilon,\lambda)\in]0,1[\times]0,150[$.}
	}
		\label{fig1}
\end{figure}

\begin{figure}[ht!]
  \centering
   % \reflectbox{%
      \includegraphics[width=0.6\textwidth]{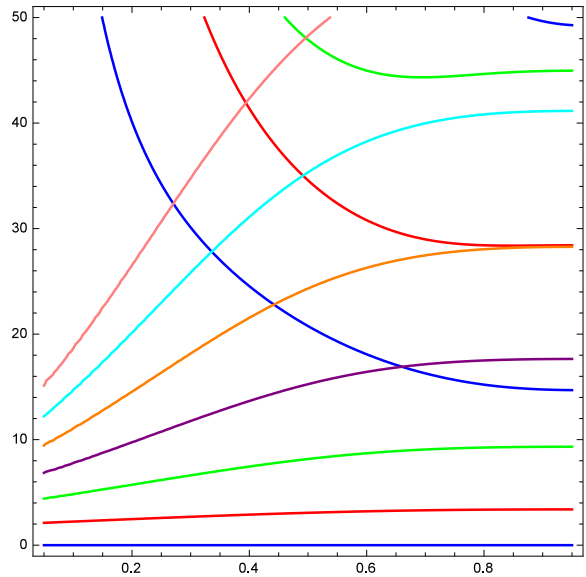}
			%}
  \caption{{Solution branches of equation \eqref{implicitformula1N} with $N=2$, $M=\pi$ for $(\varepsilon,\lambda)\in]0,1[\times]0,50[$ . The colors refer to the choice of $l$ in \eqref{implicitformula1N}: blue ($l=0$), red ($l=1$), green ($l=2$), purple ($l=3$), orange ($l=4$), cyan ($l=5$), pink ($l=6$).}}

	%Branches of eigenvalues in the case $N=2$, $M=\pi$. The colors refer to the index $l\in\mathbb N$, in particular blue ($l=0$), red ($l=1$), green ($l=2$), purple ($l=3$), orange ($l=4$), cyan ($l=5$), pink ($l=6$) . The range considered is $(\varepsilon,\lambda)\in]0,1[\times]0,50[$.}}
          % looking the other way!}
					\label{fig2}
\end{figure}

It is interesting to compare our results with those in  \cite{niwa}, where authors consider the Neumann Laplacian  in the annulus $1-\varepsilon <|x|<1$ and prove that for $N=2$ the first positive eigenvalue  is a decreasing function of $\varepsilon$. 
We note that our analysis concerns all eigenvalues $\lambda_l$ with arbitrary indexes and multiplicity, and  that  we do not prove global monotonocity 
 of $\lambda_l(\varepsilon)$, which in fact does not hold  for any $l$; see Figures \ref{fig1}, \ref{fig2}. 

The proof of our results relies on the use of Bessel functions which allows to recast  problem (\ref{Neu}) in the form of an equation  $F(\lambda , \varepsilon)=0$ in the unknowns $\lambda ,\varepsilon$.  Then, after some preparatory work, it is possible to apply  the Implicit Function Theorem  and conclude. We note that, despite the idea of the proof is rather simple and used also in other contexts (see e.g., \cite{lape}), the rigorous application of this  method requires  lenghty computations, suitable Taylor's expansions and estimates for the corresponding remainders, as well as recursive formulas for the cross-products of Bessel functions and their derivatives.

Importantly, the multiplicity of the eigenvalues which is often an obstruction in the application of standard asymptotic analysis, does not affect our method. 

We note that if the ball $B$ is replaced by a general bounded smooth domain $\Omega$, the convergence of the Neumann eigenvalues to the Steklov eigenvalues  when the mass
concentrates in a neighborhood of $\partial \Omega $ still holds. However, the explicit computation of the appropriate  formula generalizing (\ref{intro2}) is not easy and requires a 
completely different technique which will be discussed in a forthcoming paper. 

We also  note that an asymptotic analysis of similar but different problems is contained in \cite{naza1, naza2}, where  by the way explicit computations of the coefficients  in the asymptotic expansions of the eigenvalues are not provided.  

It would be interesting to investigate the monotonicity properties of the Neumann eigenvalues in the case of more general families of mass densities $\rho_{\varepsilon}$.  However, we believe that  it would be  difficult to adapt our method (which is based on explicit representation formulas) even in the case of radial mass densities (note that if $\rho_{\varepsilon}$ is not radial one could obtain a limiting Steklov-type problem with non-constant mass density, see \cite{arr} for a general discussion).  

This paper is organized as follows. The proof of formula (\ref{intro2}) is discussed in Section \ref{sec:2}.  In particular, Subsection \ref{sec:3} is devoted to   certain technical estimates which are necessary for the rigorous justification of our arguments. In Subsection \ref{sec:4} we consider also the case $N=1$ and prove formula (\ref{intro2}) for $\lambda_1$ which, by the way, is the only non zero eigenvalue of the one dimensional Steklov problem.
In  Appendix we establish the required recursive formulas for the cross-products of Bessel functions and their derivatives which are deduced by the standard 
formulas available in the literature.

%%%%%%%%%%%%%%%%%%%%%%%%%%%%%%%%%%%%%%%%%%%%%%%%%%%%%%%%%%%%%%%%%%%%%%%%%%%%%%%%%%%%%%%%%%%%%%%%%%%%
%%%%%%%%%% ASYMPTOTICS OF THE EIGENVALUES
%%%%%%%%%%%%%%%%%%%%%%%%%%%%%%%%%%%%%%%%%%%%%%%%%%%%%%%%%%%%%%%%%%%%%%%%%%%%%%%%%%%%%%%%%%%%%%%%%%%%
%%%%%%%%%%

\section{Asymptotic behavior of Neumann eigenvalues}\label{sec:2}

It is convenient to use the standard spherical coordinates $(r,\theta)$ in ${\mathbb{R}}^N$, where $\theta=(\theta_1,...\theta_{N-1})$. The corresponding trasformation of coordinates is
\begin{eqnarray*}
x_1&=&r \cos(\theta_1),\\
x_2&=&r\sin(\theta_1)\cos(\theta_2),\\
\vdots\\
x_{N-1}&=&r\sin(\theta_1)\sin(\theta_2)\cdots\sin(\theta_{N-2})\cos(\theta_{N-1}),\\
x_N&=&r\sin(\theta_1)\sin(\theta_2)\cdots\sin(\theta_{N-2})\sin(\theta_{N-1}),
\end{eqnarray*}
with $\theta_1,...,\theta_{N-2}\in [0,\pi]$, $\theta_{N-1}\in [0,2\pi[$ (here it is understood that $\theta_1\in [0,2\pi[$ if $N=2$). We denote by $\delta$ the Laplace-Beltrami operator on the unit sphere $\mathbb S^{N-1}$ of $\mathbb R^N$, which can be written in spherical coordinates as 
\begin{equation*}
\delta=\sum_{j=1}^{N-1}\frac{1}{q_j(\sin{\theta_j})^{N-j-1}}\frac{\partial}{\partial\theta_j}\left((\sin{\theta_j})^{N-j-1}\frac{\partial}{\partial\theta_j}\right),
\end{equation*}
where
$$
q_1=1,\ \ \ \ q_j=(\sin{\theta_1}\sin{\theta_2}\cdots\sin_{\theta_{j-1}})^2,\ \ \ \ j=2,...,N-1,
$$
see e.g., \cite[p. 40]{koz}.
 To shorten notation, in what follows  we will denote by $a$ and $b$ the quantities defined by 
$$a=\sqrt{\lambda\varepsilon}(1-\varepsilon),\ \ {\rm and}\ \  b=\sqrt{\lambda\tilde\rho_{\varepsilon}}(1-\varepsilon),$$ where
\begin{equation*}
\tilde\rho_{\varepsilon}=\frac{M-\varepsilon{\omega_N}\left(1-\varepsilon\right)^N}{\omega_N\left(1-\left(1-\varepsilon\right)^N\right)}.
\end{equation*}
As customary, we denote by $J_{\nu}$ and $Y_{\nu}$ the Bessel functions of the first and second species and order $\nu$ respectively (recall that $J_{\nu}$ and $Y_{\nu}$ are solutions of the Bessel equation
$
z^2y''(z)+zy'(z)+(z^2-\nu^2)y(z)=0
$).

We begin with the  following lemma.

\begin{lem}\label{solutionsN}
Given an eigenvalue $\lambda$ of problem (\ref{Neu}), a corresponding eigenfunction  $u$  is of the form $u(r,\theta)=S_l(r)H_l(\theta)$ where $H_l(\theta)$ is a spherical harmonic of some order $l\in\mathbb N$ and
\small
\begin{equation}\label{R}
S_l(r)=\left\{
\begin{array}{ll}
r^{1-\frac{N}{2}}J_{\nu_l}(\sqrt{\lambda\varepsilon} r),& {\rm if\ }r<1-\varepsilon ,  \\\ \\
r^{1-\frac{N}{2}}\left(\alpha J_{\nu_l}(\sqrt{\lambda\tilde\rho_{\varepsilon}} r)+\beta Y_{\nu_l}(\sqrt{\lambda\tilde\rho_{\varepsilon}} r)\right),& {\rm if\ }1-\varepsilon<r<1,
\end{array}
\right.
\end{equation}
\normalsize
where $\nu_l=\frac{(N+2l-2)}{2}$ and $\alpha$, $\beta$ are given by
\begin{eqnarray*}
{\alpha}={\frac{\pi b}{2}}\left({J_{\nu_l}(a)Y_{\nu_l}'(b)-\frac{a}{b}J_{\nu_l}'(a)Y_{\nu_l}(b)}\right) ,\\
{\beta}={\frac{\pi b}{2}}\left({{\frac{a}{b}}J_{\nu_l}(b)J_{\nu_l}'(a)-J_{\nu_l}'(b)J_{\nu_l}(a)}\right).\nonumber
\end{eqnarray*}
\end{lem}
\proof
Recall that the Laplace operator can be written in spherical coordinates as 
$$\Delta=\partial_{rr}+\frac{N-1}{r}\partial_r+\frac{1}{r^2}\delta.$$
In order to solve the equation $-\Delta u= \lambda \rho_{\varepsilon}u$, we  separate variables so that $u(r,\theta)=S(r)H(\theta)$.  Then using $l(l+N-2)$, $l\in\mathbb N$, as separation constant, we obtain  the equations
\begin{equation}\label{radial}
r^2S''+r(N-1)S'+r^2\lambda\rho_{\varepsilon}S-l(l+N-2)S=0
\end{equation}
and
\begin{equation}\label{angular}
-\delta H=l(l+N-2)H.
\end{equation}
By setting $S(r)=r^{1-\frac{N}{2}}\tilde S(r)$ into (\ref{radial}), it follows that $\tilde S(r)$ satisfies the Bessel equation
\begin{equation*}
\tilde S''+\frac{1}{r}\tilde S'+\left(\lambda\rho_{\varepsilon}-\frac{\nu_l^2}{r^2}\right)\tilde S=0.
\end{equation*}
Since solutions $u$ of (\ref{Neu}) are bounded on $\Omega$  and $Y_{\nu_l}(z)$ blows up at $z=0$, it follows that for $r<1-\varepsilon$, $S(r)$ is a multiple of the function $r^{1-\frac{N}{2}}J_{\nu_l}(\sqrt{\lambda\varepsilon}r)$. For $1-\varepsilon<r<1$, $S(r)$ is  a linear combination of the functions $r^{1-\frac{N}{2}}J_{\nu_l}(\sqrt{\lambda\tilde{\rho_{\varepsilon}}}r)$ and $ r^{1-\frac{N}{2}}Y_{\nu_l}(\sqrt{\lambda\tilde{\rho_{\varepsilon}}}r)$. On the other hand, the solutions of (\ref{angular}) are the spherical harmonics of order $l$. Then $u$ can be written as in (\ref{R}) for suitable values of $\alpha , \beta\in {\mathbb{R}}$.

Now we compute the coefficients $\alpha$ and $\beta$ in (\ref{R}). Since the right-hand side  of the equation in (\ref{Neu}) is a function in $L^2(\Omega)$ then by standard regularity theory a solution $u$ of (\ref{Neu}) belongs to the standard Sobolev space $H^2(\Omega)$, hence $\alpha$ and $\beta$ must be chosen in such a way that $u$ and $\partial_r u$ are continuous at $r=1-\varepsilon$, that is
\begin{equation*}\left\{
\begin{array}{ll}
\alpha J_{\nu_l}(\sqrt{\lambda\tilde\rho_{\varepsilon}}(1-\varepsilon))+\beta Y_{\nu_l}(\sqrt{\lambda\tilde\rho_{\varepsilon}}(1-\varepsilon))=J_{\nu_l}(\sqrt{\lambda\varepsilon}(1-\varepsilon))\,,\\
\alpha J_{\nu_l}'(\sqrt{\lambda\tilde\rho_{\varepsilon}}(1-\varepsilon))+\beta Y_{\nu_l}'(\sqrt{\lambda\tilde\rho_{\varepsilon}}(1-\varepsilon))=\sqrt{\frac{\varepsilon}{\tilde\rho_{\varepsilon}}}J_{\nu_l}'(\sqrt{\lambda\varepsilon}(1-\varepsilon))\,.
\end{array}
\right.
\end{equation*}
Solving the system we obtain

\begin{eqnarray*}
\alpha=\frac{J_{\nu_l}(a)Y_{\nu_l}'(b)-\frac{a}{b}J_{\nu_l}'(a)Y_{\nu_l}(b)}{J_{\nu_l}(b)Y_{\nu_l}'(b)-J_{\nu_l}'(b)Y_{\nu_l}(b)}\ ,\ \ \ \ \beta=\frac{\frac{a}{b}J_{\nu_l}(b)J_{\nu_l}'(a)-J_{\nu_l}'(b)J_{\nu_l}(a)}{J_{\nu_l}(b)Y_{\nu_l}'(b)-J_{\nu_l}'(b)Y_{\nu_l}(b)}.
\end{eqnarray*}
Note that $J_{\nu_l}(b)Y_{\nu_l}'(b)-J_{\nu_l}'(b)Y_{\nu_l}(b)$ is the Wronskian in $b$, which is known to be $\frac{2}{\pi b}$ (see \cite[\S 9]{abram}). This concludes the proof.
\endproof

We are ready to establish an implicit characterization of the eigenvalues of (\ref{Neu}).

\begin{prop}\label{implicit1N}
The nonzero eigenvalues $\lambda$ of problem (\ref{Neu}) are given implicitly as zeros of  the equation 
\begin{equation}\label{implicitformula1N}
\left(1-\frac{N}{2}\right)P_1(a,b)+\frac{b}{(1-\varepsilon)}P_2(a,b)=0
\end{equation}
where

\normalsize
\begin{eqnarray*}
P_1(a,b)&=&J_{\nu_l}(a)\left(Y_{\nu_l}'(b)J_{\nu_l}(\frac{b}{1{-}\varepsilon}){-}J_{\nu_l}'(b)Y_{\nu_l}(\frac{b}{1{-}\varepsilon})\right)\\
&{+}&\frac{a}{b}J_{\nu_l}'(a)\left(J_{\nu_l}(b)Y_{\nu_l}(\frac{b}{1{-}\varepsilon}){-}Y_{\nu_l}(b)J_{\nu_l}(\frac{b}{1{-}\varepsilon})\right),\\
P_2(a,b)&=&J_{\nu_l}(a)\left(Y_{\nu_l}'(b)J_{\nu_l}'(\frac{b}{1{-}\varepsilon}){-}J_{\nu_l}'(b)Y_{\nu_l}'(\frac{b}{1{-}\varepsilon})\right)\\
&{+}&\frac{a}{b}J_{\nu_l}'(a)\left(J_{\nu_l}(b)Y_{\nu_l}'(\frac{b}{1{-}\varepsilon}){-}Y_{\nu_l}(b)J_{\nu_l}'(\frac{b}{1{-}\varepsilon})\right).
\end{eqnarray*}
\normalsize

\proof 
By Lemma \ref{solutionsN}, an eigenfunction $u$ associated with an eigenvalue $\lambda $ is of the form $u(r,\theta )=S_l(r)H_l(\theta)$ where  for $r>1-\varepsilon$
\normalsize
\begin{eqnarray*}
S_l(r)&{=}&\frac{\pi b}{2}r^{1{-}\frac{N}{2}}\left[\left(J_{\nu_l}(a)Y_{\nu_l}'(b){-}\frac{a}{b}J_{\nu_l}'(a)Y_{\nu_l}(b)\right)J_{\nu_l}(\frac{b r}{1-\varepsilon})\right.\\
&{+}&\left.\left(\frac{a}{b}J_{\nu_l}(b)J_{\nu_l}'(a){-}J_{\nu_l}'(b)J_{\nu_l}(a)\right)Y_{\nu_l}(\frac{b r}{1-\varepsilon})\right].
\end{eqnarray*}
\normalsize
We require that $\frac{\partial u}{\partial\nu}=\frac{\partial u}{\partial r}_{|_{r=1}}=0$, which is true if and only if
\normalsize													
\begin{multline*}
\frac{\pi b}{2}\left(1{-}\frac{N}{2}\right)\left[\left(J_{\nu_l}(a)Y_{\nu_l}'(b){-}\frac{a}{b}J_{\nu_l}'(a)Y_{\nu_l}(b)\right)J_{\nu_l}(\frac{b}{1-\varepsilon})\right.\\
{+}\left.\left(\frac{a}{b}J_{\nu_l}(b)J_{\nu_l}'(a){-}J_{\nu_l}'(b)J_{\nu_l}(a)\right)Y_{\nu_l}(\frac{b}{1-\varepsilon})\right]\\
{+}\frac{\pi b^2}{2(1-\varepsilon)}\left[\left(J_{\nu_l}(a)Y_{\nu_l}'(b){-}\frac{a}{b}J_{\nu_l}'(a)Y_{\nu_l}(b)\right)J_{\nu_l}'(\frac{b}{1-\varepsilon})\right.\\
{+}\left.\left(\frac{a}{b}J_{\nu_l}(b)J_{\nu_l}'(a){-}J_{\nu_l}'(b)J_{\nu_l}(a)\right)Y_{\nu_l}'(\frac{b}{1-\varepsilon})\right]=0.
\end{multline*}
\normalsize
The previous equation can be clearly rewritten in the form (\ref{implicitformula1N}).
\endproof
\end{prop}

We now prove the following.

\begin{lem}\label{implicitepsN}
Equation (\ref{implicitformula1N}) can be written in the form
\begin{eqnarray}\label{simplified}\lefteqn{
\lambda^2\varepsilon\left(\frac{M}{3N\omega_N}-\frac{1}{\nu_l(1+\nu_l)}\right)     
+\lambda\varepsilon \left( \frac{N}{2}-\nu_l+\frac{(2-N)N\omega_N}{2\nu_l(1+\nu_l)M}  \right)   -2\lambda + 
\frac{2N\omega_Nl}{M} }\nonumber \\
& &\qquad\qquad\qquad\qquad\qquad\qquad\quad-\frac{2N\omega_Nl}{M}\left(\frac{N-1}{2}-\frac{\omega_N}{M} -\nu_l \right)\varepsilon
+ \mathcal R(\lambda,\varepsilon)=0
\end{eqnarray}
where  $\mathcal R(\lambda , \varepsilon )=O(\varepsilon \sqrt {\varepsilon})$ as $\varepsilon\to 0$.
\proof

We plan to divide the left-hand side of (\ref{implicitformula1N}) by $J_{\nu_l}'(a)$ and to analyze the resulting terms using the known Taylor's series for Bessel functions. 
Note that  $J_{\nu_l}'(a)>0$ for all $\varepsilon$ small enough. We split our analysis into three steps. 

{\it Step 1.} We consider the term $\frac{P_2(a,b)}{J_{\nu_l}'(a)}$, that is
\small
\begin{multline}\label{implicit2}
\frac{J_{\nu_l}(a)}{J_{\nu_l}'(a)}\left[Y_{\nu_l}'(b)J_{\nu_l}'(\frac{b}{1-\varepsilon})-Y_{\nu_l}'(\frac{b}{1-\varepsilon})J_{\nu_l}'(b)\right]\\
+\frac{a}{b}\left[Y_{\nu_l}'(\frac{b}{1-\varepsilon})J_{\nu_l}(b)-Y_{\nu_l}(b)J_{\nu_l}'(\frac{b}{1-\varepsilon})\right].
\end{multline}
\normalsize

Using Taylor's formula, we write the derivatives of the Bessel functions in (\ref{implicit2}), call them ${\mathcal{C}}'_{\nu_l}$, as follows
\begin{equation}\label{tayc}
{\mathcal{C}}'_{\nu_l}\left(\frac{b}{1-\varepsilon} \right)={\mathcal{C}}'_{\nu_l}(b)+{\mathcal{C}}''_{\nu_l}(b)\frac{\varepsilon b}{1-\varepsilon}+\dots +\frac{{\mathcal{C}}^{(n)}_{\nu_l}(b)}{(n-1)!}\left(\frac{\varepsilon b}{1-\varepsilon}\right)^{n-1}+o\left(\frac{\varepsilon b}{1-\varepsilon}\right)^{n-1}. 
\end{equation}
Then, using (\ref{tayc}) with $n=4$ for $J_{\nu_l}'$ and $Y_{\nu_l}'$ we get 
\small
\begin{multline}\label{numero1}
\frac{J_{\nu_l}(a)}{J_{\nu_l}'(a)}\left[\frac{\varepsilon b}{1-\varepsilon}\left(Y_{\nu_l}'(b)J_{\nu_l}''(b)-J_{\nu_l}'(b)Y_{\nu_l}''(b)\right)\right.
+\frac{\varepsilon^2 b^2}{2(1-\varepsilon)^2}\left(Y_{\nu_l}'(b)J_{\nu_l}'''(b)-J_{\nu_l}'(b)Y_{\nu_l}'''(b)\right)\\
\left.+\frac{\varepsilon^3 b^3}{6(1-\varepsilon)^3}\left(Y_{\nu_l}'(b)J_{\nu_l}''''(b)-J_{\nu_l}'(b)Y_{\nu_l}''''(b)\right)+R_1(b)\right]\\
+\frac{a}{b}\left[\left(J_{\nu_l}(b)Y_{\nu_l}'(b)-Y_{\nu_l}(b)J_{\nu_l}'(b)\right)+\frac{\varepsilon b}{1-\varepsilon}\left(J_{\nu_l}(b)Y_{\nu_l}''(b)-Y_{\nu_l}(b)J_{\nu_l}''(b)\right)\right.\\
\left.+\frac{\varepsilon^2 b^2}{2(1-\varepsilon)^2}\left(J_{\nu_l}(b)Y_{\nu_l}'''(b)-Y_{\nu_l}(b)J_{\nu_l}'''(b)\right)+R_2(b)\right],
\end{multline}
\normalsize
where 
\begin{equation}
\label{erre1}
R_1(b)=\sum_{k=4}^{+\infty}\frac{\varepsilon^kb^k}{k!(1-\varepsilon)^k}\biggl(Y'_{\nu_{l}}(b)J^{(k+1)}_{\nu_l}(b)-J'_{\nu_l}(b) Y_{\nu_{l}}^{(k+1)} (b) \biggr)\, 
\end{equation}
and
\begin{equation}
\label{erre2}
R_2(b)=\sum_{k=3}^{+\infty}\frac{\varepsilon^kb^k}{k!(1-\varepsilon)^k}\biggl(J_{\nu_{l}}(b)Y^{(k+1)}_{\nu_l}(b)-Y_{\nu_l}(b) J_{\nu_{l}}^{(k+1)} (b)\biggr)\, . 
\end{equation}
Let  $ R_3$ be the remainder defined in  Lemma \ref{Jalemma}. We set 
\small
\begin{multline}\label{numero2}
R(\lambda,\varepsilon)=R_3(a)\left[\frac{\varepsilon b}{1-\varepsilon}\left(Y_{\nu_l}'(b)J_{\nu_l}''(b){-}J_{\nu_l}'(b)Y_{\nu_l}''(b)\right)\right.\\
{+}\left.\frac{\varepsilon^2 b^2}{2(1-\varepsilon)^2}\left(Y_{\nu_l}'(b)J_{\nu_l}'''(b){-}J_{\nu_l}'(b)Y_{\nu_l}'''(b)\right)\right.\\
\left.{+}\frac{\varepsilon^3 b^3}{6(1-\varepsilon)^3}\left(Y_{\nu_l}'(b)J_{\nu_l}''''(b){-}J_{\nu_l}'(b)Y_{\nu_l}''''(b)\right)\right]\\
{+}R_1(b)\left[\frac{a}{{\nu_l}}{+}\frac{a^3}{2{\nu_l}^2(1+{\nu_l})}\right]
{+}R_2(b)\frac{a}{b}{+}R_3(a)R_1(b).
\end{multline}
\normalsize
By Lemma \ref{remainder2}, it turns out that $R(\lambda,\varepsilon)= O(\varepsilon^3)$ as $\varepsilon\rightarrow 0$. 

We also set
\begin{eqnarray*}
&&f(\varepsilon)=b_1^2(\varepsilon)a_1^3(\varepsilon)f_1(\varepsilon);\\
&&g(\varepsilon)=b_1^2(\varepsilon)a_1(\varepsilon)g_1(\varepsilon)+a_1^3(\varepsilon)g_2(\varepsilon);\\
&&h(\varepsilon)=a_1(\varepsilon)h_1(\varepsilon)+\varepsilon^2\frac{a_1^3(\varepsilon)}{b_1^2(\varepsilon)}h_2(\varepsilon);\\
&&k(\varepsilon)=\frac{a_1(\varepsilon)}{b_1^2(\varepsilon)}k_1(\varepsilon),
\end{eqnarray*}
where 
\small
\begin{eqnarray*}
&&a_1(\varepsilon)=\frac{a}{\sqrt{\lambda \varepsilon}} =(1-\varepsilon);\\
&&b_1(\varepsilon)= b\sqrt{\frac{\varepsilon }{\lambda }} ;\\
&&f_1(\varepsilon)=\frac{1}{6{\nu_l}^2(1+{\nu_l})(1-\varepsilon)^3};\\
&&g_1(\varepsilon)=\frac{1}{3{\nu_l}(1-\varepsilon)^3};\\
&&g_2(\varepsilon)=-\frac{1}{{\nu_l}^2(1+{\nu_l})(1-\varepsilon)}+\frac{\varepsilon}{2{\nu_l}^2(1+{\nu_l})(1-\varepsilon)^2}-\frac{\varepsilon^2(3+2{\nu_l}^2)}{6{\nu_l}^2(1+{\nu_l})(1-\varepsilon)^3};\\
&&h_1(\varepsilon)=-\frac{2}{{\nu_l}(1-\varepsilon)}+\frac{\varepsilon}{{\nu_l}(1-\varepsilon)^2}-\frac{\varepsilon^2(3+2{\nu_l}^2)}{3{\nu_l}(1-\varepsilon)^3}-\frac{\varepsilon}{(1-\varepsilon)^2};\\
&&h_2(\varepsilon)=\frac{1}{(1+{\nu_l})(1-\varepsilon)}-\frac{3\varepsilon}{2(1+{\nu_l})(1-\varepsilon)^2}+\frac{\varepsilon^2({\nu_l}^4+11{\nu_l}^2)}{6{\nu_l}^2(1+{\nu_l})(1-\varepsilon)^3};\\
&&k_1(\varepsilon)=2+\frac{2\varepsilon {\nu_l}}{(1-\varepsilon)}-\frac{3\varepsilon^2 {\nu_l}}{(1-\varepsilon)^2}+\frac{\varepsilon^3({\nu_l}^4+11{\nu_l}^2)}{3 \nu_l(1-\varepsilon)^3}-\frac{2\varepsilon}{(1-\varepsilon)}+\frac{\varepsilon^2(2+{\nu_l}^2)}{(1-\varepsilon)^2}.\\
\end{eqnarray*}
\normalsize
Note that functions $f,g,h,k$ are continuous at $\varepsilon =0$ and $f(0), g(0), h(0), k(0)\ne 0$.

Using the explicit formulas for the cross products of Bessel functions given by Lemma \ref{inducross} and Corollary \ref{formulaused} in (\ref{numero1}),  (\ref{implicit2}) can be written as
\begin{equation}\label{implicit5}
\frac{1}{\sqrt{\lambda}\pi}\varepsilon\sqrt{\varepsilon}k(\varepsilon)+\frac{\sqrt{\lambda}}{\pi}\varepsilon\sqrt{\varepsilon}h(\varepsilon)+\frac{\lambda\sqrt{\lambda}}{\pi}\varepsilon^2\sqrt{\varepsilon}g(\varepsilon)+\frac{\lambda^2\sqrt{\lambda}}{\pi}\varepsilon^3\sqrt{\varepsilon}f(\varepsilon)+R(\lambda,\varepsilon).
\end{equation}

{\it Step 2.} We consider  the quantity $\frac{P_1(a,b)}{J_{\nu_l}'(a)}$, that is 
\small
\begin{multline}\label{firstsum}
\frac{J_{\nu_l}(a)}{J_{\nu_l}'(a)}\left[Y_{\nu_l}'(b)J_{\nu_l}(\frac{b}{1-\varepsilon})-J_{\nu_l}'(b)Y_{\nu_l}(\frac{b}{1-\varepsilon})\right]\\
+\frac{a}{b}\left[J_{\nu_l}(b)Y_{\nu_l}(\frac{b}{1-\varepsilon})-Y_{\nu_l}(b)J_{\nu_l}(\frac{b}{1-\varepsilon})\right].
\end{multline}
\normalsize
Proceeding as in Step 1 and   setting

\small
\begin{eqnarray*}
&&\tilde f(\varepsilon)=-\frac{a_1^3(\varepsilon)b_1(\varepsilon)}{2\pi{\nu_l}^2(1+{\nu_l})(1-\varepsilon)^2};\\
&&\tilde g(\varepsilon)=\frac{a_1^3(\varepsilon)}{b_1(\varepsilon)}\left(\frac{1}{\pi{\nu_l}^2(1+{\nu_l})}+\frac{\varepsilon^2}{2\pi(1+{\nu_l})(1-\varepsilon)^2}\right)-\frac{a_1(\varepsilon)b_1(\varepsilon)}{{\nu_l}\pi(1-\varepsilon)^2};\\
&&\tilde h(\varepsilon)=\frac{a_1(\varepsilon)}{b_1(\varepsilon)}\left(\frac{2}{{\nu_l}\pi}+\frac{2\varepsilon}{\pi(1-\varepsilon)}+\frac{({\nu_l}-1)}{\pi(1-\varepsilon)^2}\varepsilon^2\right),\\
\end{eqnarray*}
\normalsize
one can prove that  (\ref{firstsum}) can be written as
\begin{equation}\label{part1}
\varepsilon \tilde h(\varepsilon)+\lambda\varepsilon^2 \tilde g(\varepsilon)+\lambda^2\varepsilon^3 \tilde f(\varepsilon)+\hat R(\lambda,\varepsilon),
\end{equation}
where   $\hat R({\lambda,\varepsilon})= O(\varepsilon^2\sqrt{\varepsilon})$ as $\varepsilon\rightarrow 0$; see  Lemma \ref{remainder2}. 

{\it Step 3.}
 We combine (\ref{implicit5}) and (\ref{part1}) and rewrite equation (\ref{implicitformula1N}) in the form
\small
\begin{multline}
\varepsilon(1-\frac{N}{2})\tilde h(\varepsilon)+\varepsilon\frac{ b_1(\varepsilon)k(\varepsilon)}{\pi (1-\varepsilon)}+\lambda\varepsilon^2(1-\frac{N}{2}) \tilde g(\varepsilon)+\lambda\varepsilon\frac{b_1(\varepsilon)h(\varepsilon)}{\pi (1-\varepsilon)}\label{implicitN}\\
+\lambda^2\varepsilon^3(1-\frac{N}{2})\tilde f(\varepsilon)+\lambda^2\varepsilon^2\frac{b_1(\varepsilon)g(\varepsilon)}{\pi (1-\varepsilon)}+\lambda^3\varepsilon^3\frac{b_1(\varepsilon)f(\varepsilon)}{\pi (1-\varepsilon)}+\mathcal R_0(\lambda,\varepsilon)=0,
\end{multline}
\normalsize
where
\begin{eqnarray*}
&&\mathcal R_0(\lambda,\varepsilon)=\frac{\sqrt{\lambda}b_1(\varepsilon)}{(1-\varepsilon)\sqrt{\varepsilon}}R(\lambda,\varepsilon)+\left(1-\frac{N}{2}\right)\hat R(\lambda,\varepsilon).
\end{eqnarray*}
Note that $\mathcal R_0(\lambda,\varepsilon) =O(\varepsilon^2\sqrt{\varepsilon})$ as $\varepsilon\rightarrow 0$. Dividing by $\varepsilon $ in (\ref{implicitN}) and setting $\mathcal R_1(\lambda,\varepsilon)=\frac{\mathcal R_0(\lambda,\varepsilon)}{\varepsilon}$, we obtain 
\begin{eqnarray}
&&(1-\frac{N}{2})\tilde h(\varepsilon)+\frac{ b_1(\varepsilon)k(\varepsilon)}{\pi (1-\varepsilon)}+\lambda\varepsilon(1-\frac{N}{2}) \tilde g(\varepsilon)+\lambda\frac{b_1(\varepsilon)h(\varepsilon)}{\pi (1-\varepsilon)}\label{implicitN1}\\
&&+\lambda^2\varepsilon^2(1-\frac{N}{2})\tilde f(\varepsilon)+\lambda^2\varepsilon\frac{b_1(\varepsilon)g(\varepsilon)}{\pi (1-\varepsilon)}+\lambda^3\varepsilon^2\frac{b_1(\varepsilon)f(\varepsilon)}{\pi (1-\varepsilon)}+\mathcal R_1(\lambda,\varepsilon)=0\nonumber .
\end{eqnarray}

We now multiply in \eqref{implicitN1} by $\frac{\pi\nu_l (1-\varepsilon) }{ b_1(\varepsilon)}$ which is a positive quantity  for all $0<\varepsilon<1$. Taking into account the definitions of functions $g,h,k, \tilde g,\tilde h$,  we can finally rewrite (\ref{implicitN1}) in the form
\begin{eqnarray}\lefteqn{
\lambda^2\varepsilon\left(\frac{\hat\rho (\varepsilon )}{3}-\frac{1}{\nu_l(1+\nu_l)}\right)     
+\lambda\varepsilon \left( \frac{N}{2}-\nu_l+\frac{2-N}{2\nu_l(1+\nu_l)\hat\rho (\varepsilon )}  \right)   -2\lambda  }\nonumber \\
& &\qquad\qquad\qquad\qquad\qquad\qquad\qquad\qquad\qquad\qquad+ \frac{2l\left(1+\varepsilon\nu_l\right)}{\hat\rho (\varepsilon )} + \mathcal R(\lambda,\varepsilon)=0 ,
\end{eqnarray}
where 
$$\hat \rho (\varepsilon )=\varepsilon \tilde\rho (\varepsilon )=\frac{M-\omega_N\varepsilon(1-\varepsilon)^N}{\omega_N\left(N-\frac{N(N-1)}{2}\varepsilon-\sum_{k=3}^N\binom{N}{k}(-1)^k\varepsilon^{k-1}\right)},
$$
and $\mathcal R(\lambda,\varepsilon) =O(\varepsilon\sqrt {\varepsilon} )$ as $\varepsilon \to 0$.
The formulation in (\ref{simplified})  can be easily deduced by observing that 
$$
\hat \rho_{\varepsilon}= \frac{M}{N\omega_N}+2\frac{M}{N\omega_N}\left( \frac{N-1}{4}-\frac{\omega_N}{2M}  \right)\varepsilon +O(\varepsilon^2),\ \ {\rm as}\ \varepsilon \to 0 .
$$

\endproof 
\end{lem}

We are now ready to prove our main result

\begin{thm}\label{derN}
All  eigenvalues of problem (\ref{Neu}) have the following asymptotic behavior
\begin{equation}\label{asymptotic}
\lambda_l(\varepsilon)= \lambda_l+\left(\frac{2l\lambda_l}{3}+\frac{2\lambda^2_l}{N(2l+N)}\right)\varepsilon+o(\varepsilon ),\ \ \ {\rm as}\ \varepsilon \to 0,
\end{equation}
where $\lambda_l$ are the eigenvalues of problem (\ref{Ste}).

Moreover, for each $l\in {\mathbb{N}}$ the  function defined by $\lambda_l(\varepsilon )$ for  $\varepsilon >0$ and $\lambda_l(0)=\lambda_l$, is continuous in the whole of $[0,1[$
and of class $C^1$ in a neighborhood of $\varepsilon =0$.
\proof

By using the Min-Max Principle and related standard arguments, one can easily prove that  $\lambda_l(\varepsilon)$ depends with continuity 
on $\varepsilon >0$ (cfr. \cite{laproeurasian}, see also \cite{lala2004}). Moreover, by using (\ref{intro1,5})  the maps $\varepsilon \mapsto \lambda_l(\varepsilon)$ can be extended by continuity at the point $\varepsilon =0$
by setting $\lambda_l(0)=\lambda_l$.

In order to prove differentiability of $\lambda_{l}(\varepsilon )$ around zero and the validity of (\ref{asymptotic}), we consider equation (\ref{simplified}) and apply the Implicit Function Theorem. Note that equation (\ref{simplified}) can be  written in  the form $F(\lambda , \varepsilon)=0$ where $F$ is a function of class $C^1$ in the variables $(\lambda , \varepsilon )\in ]0,\infty [\times[0,1[$, with 
\begin{eqnarray}
 F(\lambda , 0) &=& -2\lambda +\frac{2N\omega_Nl}{M},\nonumber\\
 F'_{\lambda}(\lambda ,0)& =& -2,\nonumber\\
F'_{\varepsilon}(\lambda , 0)&=&\lambda^2\left(\frac{M}{3N\omega_N}-\frac{1}{\nu_l(1+\nu_l)}\right)     
+\lambda \left( \frac{N}{2}-\nu_l+\frac{(2-N)N\omega_N}{2\nu_l(1+\nu_l)M}  \right)\nonumber\\
& &  -\frac{2N\omega_Nl}{M}\left(\frac{N-1}{2}-\frac{\omega_N}{M} -\nu_l \right)
\end{eqnarray}

By (\ref{intro1}), $\lambda_l=N\omega_Nl/M$ hence $F(\lambda_l,0)=0$. Since $F'_{\lambda}(\lambda_l ,0)\ne 0$, the Implicit Function Theorem combined with the 
continuity of the functions $\lambda_l(\cdot )$ allows to conclude that functions  $\lambda_l(\cdot )$ are of class $C^1$ around zero.  

We now compute the derivative of $\lambda_l(\cdot )$ at zero. Using the equality $N\omega_N/M=\lambda_l /l$ and recalling that $\nu_l=l+N/2-1$ we get
\begin{eqnarray}
F'_{\varepsilon}(\lambda_l , 0)&=&\lambda_l^2\left(  \frac{l}{3\lambda_l} -\frac{1}{\nu_l(1+\nu_l)} \right)+\lambda_l\left(1-l+\frac{\lambda_l (2-N)}{2l\nu_l(1+\nu_l)}   \right)
-2\lambda_l\left(\frac{1}{2}-l-\frac{\lambda_l}{Nl}  \right)\nonumber \\
&=&\lambda_l^2\left( \frac{1}{\nu_l(1+\nu_l)}\left(\frac{2-N}{2l}-1  \right)+\frac{2}{Nl}  \right)+\frac{4}{3}\lambda_l l=\frac{4\lambda^2_l}{N^2+2Nl}+ \frac{4}{3}\lambda_l l . \nonumber
\end{eqnarray}
Finally,  formula $\lambda'_l(0)=-F'_{\varepsilon}(\lambda_l , 0)/F'_{\lambda}(\lambda_l , 0)$ yields (\ref{intro2}) and the validity of (\ref{asymptotic}).
\endproof
\end{thm}

\begin{corol}\label{derminN}
For any $l\in {\mathbb{N}}\setminus \{0\}$ there exists $\delta_l$ such that the function $\lambda _l(\cdot )$ is strictly increasing in the interval $[0, \delta _l[$. In particular, $\lambda_l<\lambda_l(\varepsilon )$ for all $\varepsilon \in ]0, \delta_l[$. 
\end{corol}

%\newpage

%{\color{red} INSERT A PICTURE WITH SOME BRANCHES OF EIGENVALUES AND A REMARK WHERE WE SAY THAT GLOBAL MONOTONICITY DOES NOT HOLD IN GENERAL}

%%%%%%%%%%%%%%%%%%%%%%%%%%%%%%%%%%%%%%%%%%%%%%%%%%%%%%%%%%%%%%%%%%%%%%%%%%%%%%%%%%%%%%%%%%%%%%%%%%%%

%%%%%%%%%%%%%%%%%%%%%%%%%%%%%%%%%%%%%%%%%%%%%%%%%%%%%%%%%%%%%%%%%%%%%%%%%%%%%%%%%%%%%%%%%%%%%%%%%%%%

%%%%%%%%%%%%%%%%%%%%%%%%%%% APPENDIX - REMAINDERS OF TAYLOR EXPANSIONS AND CROSS PRODUCTS %%%%%%%%%%

%%%%%%%%%%%%%%%%%%%%%%%%%%%%%%%%%%%%%%%%%%%%%%%%%%%%%%%%%%%%%%%%%%%%%%%%%%%%%%%%%%%%%%%%%%%%%%%%%%%%

%%%%%%%%%%%%%%%%%%%%%%%%%%%%%%%%%%%%%%%%%%%%%%%%%%%%%%%%%%%%%%%%%%%%%%%%%%%%%%%%%%%%%%%%%%%%%%%%%%%%

\newpage
\subsection{Estimates for the remainders}\label{sec:3}

This subsection is devoted to the proof of a few technical estimates used in the proof of Lemma~\ref{implicitepsN}.

\begin{lem}\label{Jalemma}
The function $R_3$ defined by
\begin{equation}
\label{Jalemma0}
\frac{J_{\nu}(z)}{J_{\nu}'(z)}=\frac{z}{{\nu}}+\frac{z^3}{2{\nu}^2(1+{\nu})}+ R_3(z),
\end{equation}
 is $ O(z^5)$ as $z\rightarrow 0$.
\proof

Recall the well-known following representation of the Bessel functions of the first species

\begin{equation}\label{besselabram}
J_{\nu}(z)=\left(\frac{z}{2}\right)^{\nu}\sum_{j=0}^{+\infty}\frac{(-1)^j}{j!\Gamma(j+{\nu}+1)}\left(\frac{z}{2}\right)^{2j}.
\end{equation}

For clarity, we simply write 
\begin{equation}\label{besselabram1}
J_{\nu}(z)=z^{\nu}(a_0+a_2z^2+a_4z^4+O(z^5) ),
\end{equation}
hence 
\begin{equation}\label{besselabram2}
J_{\nu}'(z)=z^{\nu-1}(\nu a_0+(\nu +2)a_2z^2+(\nu+4 )a_4z^4+O(z^5) )
\end{equation}
where the coefficients $a_0,a_2, a_4$ are defined by (\ref{besselabram}). By (\ref{besselabram1}), (\ref{besselabram2}) and  standard computations it  follows
that
$$
\frac{J_{\nu}(z)}{J_{\nu}'(z)}=\frac{z}{\nu}-{\frac{2a_2}{\nu^2a_0}}z^3+O(z^5),
$$
which gives exactly (\ref{Jalemma0}).
\endproof
\end{lem}

%%%%%%%%%%%%%%%%%%%%%%%%%%%%%%%%%%%%%%%%%%%%%%%%%%%%%%%%%%%%%%%%%%%%%%%%%%%%%%%%%%%%%%%%%%%%%%%%%%%%

%%%%%%%%%%%%%%%%%%%%%%%%%%%%%%%%%%%%%%%%%%%%%%%%%%%%%%%%%%%%%%%%%%%%%%%%%%%%%%%%%%%%%%%%%%%%%%%%%%%%
\begin{lem}\label{remainder2} For any $\lambda >0$ the 
remainders $R(\lambda,\varepsilon)$ and $\hat R(\lambda , \varepsilon)$  defined in the proof of Lemma~\ref{implicitepsN} are  $ O(\varepsilon^3)$, $O(\varepsilon^2 \sqrt {\varepsilon })$, respectively,  as $\varepsilon\rightarrow 0$.
Moreover,  the same holds true for the corresponding partial derivatives  $\partial_{\lambda} R(\lambda,\varepsilon)$,  $\partial_{\lambda}\hat  R(\lambda,\varepsilon)$. 
\proof

First, we consider $R_3(a)=R_3(\sqrt{\lambda\varepsilon}(1-\varepsilon))$ where $R_3$ is defined in Lemma~\ref{Jalemma} and we differentiate it with respect to $\lambda$. We obtain
\small
\begin{eqnarray*}
\frac{\partial R_3(a)}{\partial\lambda }=\frac{a R_3'(a)}{2\lambda},
\end{eqnarray*}
\normalsize
hence by Lemma \ref{Jalemma} we can conclude that $R_3(a)$ and $\frac{\partial R_3(a)}{\partial\lambda }$ are $O(\varepsilon^2\sqrt{\varepsilon})$ as $\varepsilon\rightarrow 0$.

Now  consider $R_1(b)$ and $R_2(b)$ defined in (\ref{erre1}),  (\ref{erre2}).
%Suppose for the moment to know that the following equalities hold
%\small
%\begin{eqnarray*}
%Y_{\nu}(b)J_{\nu}^{(k)}(b)-J_{\nu}(b)Y_{\nu}^{(k)}(b)&=&\frac{2}{\pi b}\left({r({\nu},k)}+R_{{\nu},k}(b)\right),\\
%Y_{\nu}'(b)J_{\nu}^{(k)}(b)-J_{\nu}'(b)Y_{\nu}^{(k)}(b)&=&\frac{2}{\pi b}\left({q({\nu},k)}+Q_{{\nu},k}(b)\right)
%\end{eqnarray*}
%\normalsize
%for all ${\nu}\geq 0$, where $r({\nu},k)=0,1$ or $-1$, $q({\nu},k)=0,1$ or $-1$, depending on ${\nu},k$, and $Q_{{\nu},k}(b)$ and $R_{{\nu},k}(b)$ are finite sums of quotients %of the form $\frac{C({\nu},k)}{b^m}$, with $m\geq 1$ and $C({\nu},k)$ a suitable constant, depending on ${\nu},k$. This is proved in Lemma \ref{inducross}. 
 Since $\lambda>0$,  we have that $b>0$ hence  the Bessel functions are analytic in $b$ and we can write 
\small
\begin{eqnarray*}
2\sqrt{\lambda} \frac{ \partial R_1(b)}{\partial \lambda}&=&\frac{\varepsilon b_1(\varepsilon)}{\sqrt{\varepsilon}(1-\varepsilon)}\sum_{k=4}^{+\infty}\frac{b^{k-1}\varepsilon^{k-1}}{(k-1)!(1-\varepsilon)^{k-1}}\left(Y_{\nu}'(b)J_{\nu}^{(k+1)}(b)-J_{\nu}'(b)Y_{\nu}^{(k+1)}(b)\right)\\
&+&\frac{b_1(\varepsilon)}{\sqrt{\varepsilon}}\sum_{k=4}^{+\infty}\frac{\varepsilon^kb^k}{k!(1-\varepsilon)^k}\left(Y_{\nu}'(b)J_{\nu}^{(k+1)}(b)-J_{\nu}'(b)Y_{\nu}^{(k+1)}(b)\right)'.
\end{eqnarray*}
\normalsize
Here and in the sequel we  write $\nu$ instead of $\nu_l$.
Using the fact that $b=\sqrt{\lambda /\varepsilon}b_1(\varepsilon)$ and Lemma \ref{inducross} we conclude that all the cross products of the form 
$Y_{\nu}'(b)J_{\nu}^{(k+1)}(b)-J_{\nu}'(b)Y_{\nu}^{(k+1)}(b)$ and their derivatives $(Y_{\nu}'(b)J_{\nu}^{(k+1)}(b)-J_{\nu}'(b)Y_{\nu}^{(k+1)}(b))'$
are  $O(\sqrt{\varepsilon})$ and  $O(\varepsilon)$ respectively,  as $\varepsilon\rightarrow 0$.
It  follows  that $R_1(b)$ and $\partial_{\lambda}R_1(b)$ are $ O(\varepsilon^2\sqrt{\varepsilon})$ as $\varepsilon\rightarrow 0$. 

Similarly, 
\small
\begin{eqnarray*}
2\sqrt{\lambda}\frac{\partial R_2(b)}{\partial \lambda}&=&\frac{\varepsilon b_1(\varepsilon)}{\sqrt{\varepsilon}(1-\varepsilon)}\sum_{k=3}^{+\infty}\frac{b^{k-1}\varepsilon^{k-1}}{(k-1)!(1-\varepsilon)^{k-1}}\left(J_{\nu}(b)Y_{\nu}^{(k+1)}(b)-Y_{\nu}(b)J_{\nu}^{(k+1)}(b)\right)\\
&+&\frac{b_1(\varepsilon)}{\sqrt{\varepsilon}}\sum_{k=3}^{+\infty}\frac{\varepsilon^kb^k}{k!(1-\varepsilon)^k}\left(J_{\nu}(b)Y_{\nu}^{(k+1)}(b)-Y_{\nu}(b)J_{\nu}^{(k+1)}(b)\right)' ,
\end{eqnarray*}
\normalsize
hence $R_2(b)$ and $\partial_{\lambda}R_2(b)$ are $O(\varepsilon^2)$ as $\varepsilon\rightarrow 0$. 

Summing up all the terms, using Lemma~\ref{case0} and Corollary~\ref{formulaused}, we obtain 
\small
\begin{multline*}
R(\lambda,\varepsilon)=R_3(a)\left[\frac{2\varepsilon}{\pi(1-\varepsilon)}\left(\frac{{\nu}^2}{b^2}-1\right){+}\frac{\varepsilon^2}{\pi(1-\varepsilon)^2}\left(1-\frac{3{\nu}^2}{b^2}\right)\right.
\\\left.{+}\frac{\varepsilon^3b^2}{3\pi(1-\varepsilon)^3}\left(\frac{{\nu}^4+11{\nu}^2}{b^4}-\frac{3+2{\nu}^2}{b^2}+1\right)\right]\\
+R_1(b)\left[\frac{a}{{\nu}}+\frac{a^3}{2{\nu}^2(1+{\nu})}\right]
{+}R_2(b)\frac{a}{b}{+}R_3(a)R_1(b).
\end{multline*}
\normalsize
We conclude that $R(\lambda,\varepsilon)$ is $ O(\varepsilon^3)$ as $\varepsilon\rightarrow 0$. Moreover, it easily follows  that $\frac{\partial R (\lambda,\varepsilon)}{\partial \lambda}$ is also $ O(\varepsilon^3)$ as $\varepsilon\rightarrow 0$. 

The proof of the estimates for $\hat R$ and its derivatives is similar and we omit it. 
\endproof
\end{lem}

\begin{rem} According to standard Landau's notation,  saying that a function $f(z)$ is $O(g(z))$ as $z\to 0$ means that there exists $C>0$  such that 
$|f(z)|\le C|g(z)|$ for any $z$ sufficiently close to zero. Thus, using Landau's notation in the statements of Lemmas~\ref{implicitepsN}, \ref{remainder2}
understands the existence of such constants $C$, which in principle may depend on $\lambda >0$. However, a careful analysis of the proofs 
reveals that given a bounded interval of the type $[A,B]$ with $0<A<B$ then the appropriate  constants $C$ in the estimates can be taken independent
of $\lambda \in [A,B]$.
\end{rem}

\subsection{The case \texorpdfstring{$N=1$}{TEXT}}\label{sec:4}

%%%%%%%%%%%%%%%%%%%%%%%%%%%%%%    ONE DIMENSIONAL CASE %%%%%%%%%%%%%%%%%%%%%%%%%%%%%%%%%%%%%%%%%%%%%%%%%%%

%%%%%%%%%%%%%%%%%%%%%%%%%%%%%%%%%%%%%%%%%%%%%%%%%%%%%%%%%%%%%%%%%%%%%%%%%%%%%%%%%%%%%%%%%%%%%%%%%%%%%%%%%%%%

We include here a description of the case $N=1$ for the sake of completeness. Let $\Omega$ be the open interval $]-1,1[$. Problem \eqref{Ste} reads
\begin{equation}\label{Ste1}
\begin{cases}
u''(x) =0 ,& {\rm for}\  x\in]-1,1[,\\
u'(\pm 1)=\pm\lambda\frac{M}{2}u(\pm 1),
%u'(-1)=-\lambda\frac{M}{2}u(-1),\\
\end{cases}
\end{equation}
in the unknowns $\lambda$ and $u$. It is easy to see that the only eigenvalues are $\lambda_0=0$ and $\lambda_1=\frac{2}{M}$ and they are  associated  with the constant functions and the function $x$, respectively.  As in (\ref{densitaintro}), we define a mass density $\rho_{\varepsilon}$ on the whole of $]-1,1[$ by

\begin{equation*}
\rho_{\varepsilon}(x)=\left\{
\begin{array}{ll}
\frac{M}{2\varepsilon}-1+\varepsilon\,& {\rm if\ }x\in ]-1,-1+\varepsilon[\cup]1-\varepsilon,1[,\\
\varepsilon\,& {\rm if\ }x\in ]-1+\varepsilon,1-\varepsilon[.
\end{array}
\right.
\end{equation*}
Note that for any $x\in]-1,1[ $  we have $\rho_{\varepsilon}(x)\to 0$ as $\varepsilon \to 0$, and   $\int_{-1}^1\rho_{\varepsilon}dx=M$ for all $\varepsilon>0$. Problem (\ref{Neu}) for $N=1$ reads 
\begin{equation}\label{Neu1}
\left\{\begin{array}{ll}
-u''(x) =\lambda \rho_{\varepsilon}(x) u(x),\ \ & {\rm for}\  x\in]-1,1[,\\
u'(-1)=u'(1)=0. 
\end{array}\right.
\end{equation}
It is well-known from Sturm-Liouville theory that problem (\ref{Neu1}) has an increasing sequence of non-negative eigenvalues of multiplicity one. We denote the eigenvalues of \eqref{Neu1} by $\lambda_l(\varepsilon)$ with $l\in\mathbb N$.  For any $\varepsilon\in]0,1[$, the only zero eigenvalue is $\lambda_0(\varepsilon)$ and the corresponding eigenfunctions are the constant functions.

We establish an implicit characterization of the eigenvalues of \eqref{Neu1}.

\begin{prop}
The nonzero eigenvalues $\lambda$ of problem \eqref{Neu1} are given implicitly as zeros of the equation
\begin{multline}\label{implicitD1}
2\sqrt{\varepsilon\left(\frac{M}{2\varepsilon}-1+\varepsilon\right)}\cos{(2\sqrt{\lambda\varepsilon}(1-\varepsilon))}\sin{\left(2\varepsilon\sqrt{\lambda\left(\frac{M}{2\varepsilon}-1+\varepsilon\right)}\right)}\\
+\left[-\frac{M}{2\varepsilon}+1+\left(\frac{M}{2\varepsilon}-1+2\varepsilon\right)\cos{\left(2\varepsilon\sqrt{\lambda\left(\frac{M}{2\varepsilon}-1+\varepsilon\right)}\right)}\right]\sin{\left(2\sqrt{\lambda\varepsilon}(1-\varepsilon)\right)}=0.
\end{multline}
\proof
Given an eigenvalue $\lambda>0$, a solution of (\ref{Neu1}) is of the form
\small
\begin{equation*}
u(x)=\left\{
\begin{array}{ll}
A\cos{(\sqrt{\lambda\rho_2}x)}+B\sin{(\sqrt{\lambda\rho_2}x)},& {\rm for\ }x\in]-1,-1+\varepsilon[,\\\ \\
C\cos{(\sqrt{\lambda\rho_1}x)}+D\sin{(\sqrt{\lambda\rho_1}x)},& {\rm for\ }x\in]-1+\varepsilon,1-\varepsilon[,\\\ \\
E\cos{(\sqrt{\lambda\rho_2}x)}+F\sin{(\sqrt{\lambda\rho_2}x)},& {\rm for\ }x\in]1-\varepsilon,1[,
\end{array}
\right.
\end{equation*}
\normalsize
where $\rho_1=\varepsilon,\rho_2=\frac{M}{2\varepsilon}-1+\varepsilon$ and $A,B,C,D,E,F$ are suitable real numbers. We impose the continuity of $u$ and $u'$ at the points $x=-1+\varepsilon$ and $x=1-\varepsilon$ and the boundary conditions, obtaining a homogeneous system of six linear equations in six unknowns of the form $\mathcal M v=0$, where $v=(A,B,C,D,E,F)$ and $\mathcal M$ is the matrix associated with the system. We impose the condition ${\rm det}\mathcal M=0$. This yields formula \eqref{implicitD1}. 
\endproof
\end{prop}

Note that $\lambda=0$ is a solution for all $\varepsilon>0$, then we consider only the case of nonzero eigenvalues. 
Using standard Taylor's formulas, we easily prove the following

\begin{lem}
Equation \eqref{implicitD1} can be rewritten in the form
\begin{equation}\label{implicitepsD1}
M-\frac{\lambda M^2}{2}+\frac{\lambda M^2}{6}\left(1+\lambda\left(2+\frac{M}{2}\right)\right)\varepsilon+R(\lambda,\varepsilon)=0,
\end{equation}
where $R(\lambda,\varepsilon) =O(\varepsilon^2)$ as $\varepsilon\rightarrow 0$.
\end{lem}

Finally, we can prove the following theorem. 
Note that formula \eqref{formulaD1} is the same as \eqref{asymptotic} with $N=1,l=1$.

\begin{thm}
The first eigenvalue of problem \eqref{Neu1} has the following asymptotic behavior
\begin{equation}\label{formulaD1}
\lambda_1(\varepsilon)=\lambda_1+\frac{2}{3}(\lambda_1+\lambda_1^2)\varepsilon+o(\varepsilon)\ \ \ {\rm as}\ \varepsilon \to 0,
\end{equation}
where $\lambda_1=2/M$ is the only nonzero eigenvalue of problem \eqref{Ste1}. Moreover, for $l>1$ we have that $\lambda_l(\varepsilon)\rightarrow +\infty$ as $\varepsilon\rightarrow 0$.
\proof
The proof is similar to that of Theorem \ref{derN}. It is possible to prove that the eigenvalues $\lambda_l(\varepsilon)$ of \eqref{Neu1} depend with continuity on $\varepsilon>0$. We consider equation \eqref{implicitepsD1} and apply the Implicit Function Theorem. Equation \eqref{implicitepsD1} can be written in the form $F(\lambda,\varepsilon)=0$, with $F$ of class $C^1$ in $]0,+\infty[\times[0,1[$ with $F(\lambda,0)=M-\frac{\lambda M^2}{2}$, $F_{\lambda}'(\lambda,0)=-\frac{M^2}{2}$ and $F_{\varepsilon}'(\lambda,0)=\frac{\lambda M^2}{6}(1+\lambda (2+\frac{M}{2}))$.

Since  $\lambda_1=\frac{2}{M}$, $F(\lambda_1,0)= 0$ and $F_{\lambda}'(\lambda_1,0)\ne 0$, the zeros of equation (\ref{formulaD1}) in a neighborhood of $(\lambda ,0)$ are given
by the graph of a $C^1$-function $\varepsilon \mapsto \lambda (\varepsilon )$ with $\lambda (0)=\lambda_1$. We note that 
$\lambda (\varepsilon )=\lambda_1(\varepsilon )$ for all $\varepsilon $ small enough. Indeed, assuming by contradiction that $\lambda (\varepsilon )=\lambda_l(\varepsilon )$ with $l\geq 2$,  we would obtain that, possibly passing to a subsequence, $\lambda_1(\varepsilon )\to \bar \lambda $ as $\varepsilon \to 0$, for some $\bar \lambda \in [0, \lambda_1[$. Then passing to the limit in (\ref{implicitepsD1}) as $\varepsilon \to 0$ we would obtain a contradiction. Thus, $\lambda_1 (\cdot )$ is of class $C^1$ in a neighborhood of zero and $\lambda_1'(0)=-F'_{\varepsilon}(\lambda_1,0)/F'_{\lambda}(\lambda_1,0)$  which yields formula \eqref{formulaD1}. 

The divergence as $\varepsilon \to 0$ of the higher eigenvalues $\lambda_l(\varepsilon  )$ with $l>1$,
is clearly deduced by the fact that the existence of a converging subsequence of the form $\lambda_l(\varepsilon _n)$, $n\in {\mathbb{N}}$ would provide the existence of 
an eigenvalue for the limiting problem (\ref{Ste1}) different from $\lambda_0$ and $\lambda_1$, which is not admissible.

 %The fact that near $\varepsilon=0$ we say that formula \eqref{formulaD1} holds for the first eigenvalue $\lambda_1(\varepsilon)$ of problem \eqref{Neu1} follows from the fact that for each $\bar\lambda<\lambda_1$ there exists $\bar\varepsilon>0$ such that the left-hand side of \eqref{implicitepsD1} is different from zero for all $(\varepsilon,\lambda)\in]0,\bar\varepsilon[\times]0,\bar\lambda[$.
\endproof
\end{thm}

%\begin{rem}
%It is possible to show that for $l>1$, $\lambda_l(\varepsilon)\rightarrow +\infty$ as $\varepsilon\rightarrow 0$. In fact for each fixed $\bar\lambda>\lambda_1$, there exists $\bar\varepsilon>0$ such that the left-hand side of \eqref{implicitepsD1} is different from zero for all $(\varepsilon,\lambda)\in]0,\bar\varepsilon[\times]0,\bar\lambda[$.
%\end{rem}

%%%%%%%%%%%%%%%%%%%%%%%%%%%%%%%%%%%%%%%%%%%%%%%%%%%%%%%%%%%%%%%%%%%%%%%%%%%%%%%%%%%%%%%%%%%%%%%%%%%%%%%%%%%%%%%%%%%%

%%%%%%%%%%%%%%%%%%%%%%%%%%%%%%%%%%%%%%%%%%%%%%%%%%%%%%%%%%%%%%%%%%%%%%%%%%%%%%%%%%%%%%%%%%%%%%%%%%%%%%%%%%

\section{Appendix}

We provide here explicit formulas for the cross products of Bessel functions used in this paper. 

\begin{lem}\label{case0}
The following identities hold
\small
\begin{eqnarray*}
Y_{\nu}(z)J_{\nu}'(z)-J_{\nu}(z)Y_{\nu}'(z)&=&-\frac{2}{\pi z},\\
Y_{\nu}(z)J_{\nu}''(z)-J_{\nu}(z)Y_{\nu}''(z)&=&\frac{2}{\pi z^2},\\
Y_{\nu}'(z)J_{\nu}''(z)-J_{\nu}'(z)Y_{\nu}''(z)&=&\frac{2}{\pi z}\left(\frac{{\nu}^2}{z^2}-1\right),
\end{eqnarray*}
\normalsize
\proof
It is well-known (see \cite[\S 9]{abram}) that
\small
\begin{eqnarray*}J_{\nu}(z)Y_{\nu}'(z)-Y_{\nu}(z)J_{\nu}'(z)=
J_{{\nu}+1}(z)Y_{\nu}(z)-J_{\nu}(z)Y_{{\nu}+1}(z)=\frac{2}{\pi z},
\end{eqnarray*}
\normalsize
which gives  the first identity in the statement. The second identity holds since
\small
\begin{eqnarray*}
J_{\nu}(z)Y_{\nu}''(z)-Y_{\nu}(z)J_{\nu}''(z)&=&\left(J_{\nu}(z)Y_{\nu}'(z)-Y_{\nu}(z)J_{\nu}'(z)\right)'=\left(\frac{2}{\pi z}\right)'=-\frac{2}{\pi z^2}.
\end{eqnarray*}
\normalsize
The third identity holds since
\small
\begin{multline*}
Y_{\nu}'(z)J_{\nu}''(z)-J_{\nu}'(z)Y_{\nu}''(z)=Y_{\nu}'(z)\left(J_{{\nu}-1}(z)-\frac{{\nu}}{z}J_{\nu}(z)\right)'-J_{\nu}'(z)\left(Y_{{\nu}-1}(z)-\frac{{\nu}}{z}Y_{\nu}(z)\right)'\\
=Y_{\nu}'(z)J_{{\nu}-1}'(z)-J_{\nu}'(z)Y_{{\nu}-1}'(z)+\frac{{\nu}}{z^2}\left(Y_{\nu}'(z)J_{\nu}(z)-J_{\nu}'(z)Y_{\nu}(z)\right)\\
=\left(Y_{\nu}'(z)\frac{1}{2}\left(J_{{\nu}-2}(z)-J_{\nu}(z)\right)-J_{\nu}'(z)\frac{1}{2}\left(Y_{{\nu}-2}(z)-Y_{\nu}(z)\right)\right)+\frac{2{\nu}}{\pi z^3}\\
=\frac{1}{2}\left(Y_{\nu}'(z)J_{{\nu}-2}(z)-J_{\nu}'(z)Y_{{\nu}-2}(z)\right)\\
-\frac{1}{2}\left(Y_{\nu}'(z)J_{\nu}(z)-J_{\nu}'(z)Y_{\nu}(z)\right)+\frac{2{\nu}}{\pi z^3}\\
=\frac{1}{2}\left(J_{\nu}'(z)Y_{\nu}(z)-Y_{\nu}'(z)J_{\nu}(z)\right)\\
+\frac{{\nu}-1}{z}\left(Y_{\nu}'(z)J_{{\nu}-1}(z)-J_{\nu}'(z)Y_{{\nu}-1}(z)\right)-\frac{1}{\pi z}+\frac{2{\nu}}{\pi z^3}\\
=\frac{{\nu}-1}{z}\left(J_{{\nu}-1}(z)\left(Y_{{\nu}-1}(z)-\frac{{\nu}}{z}Y_{\nu}(z)\right)-Y_{{\nu}-1}(z)\left(J_{{\nu}-1}(z)-\frac{{\nu}}{z}J_{\nu}(z)\right)\right)\\
-\frac{2}{\pi z}+\frac{2{\nu}}{\pi z^3}\\
=-\frac{{\nu}({\nu}-1)}{z^2}\left(Y_{\nu}(z)J_{{\nu}-1}(z)-J_{\nu}(z)Y_{{\nu}-1}(z)\right)-\frac{2}{\pi z}+\frac{2{\nu}}{\pi z^3}\\
=\frac{2}{\pi z}\left(-1+\frac{{\nu}^2}{z^2}\right),
\end{multline*}
\normalsize
where the first, second and fourth equalities follow respectively from the well-known formulas $\mathcal C_{\nu}'(z)=\mathcal C_{{\nu}-1}(z)-\frac{{\nu}}{z}\mathcal C_{{\nu}}(z)$, $2\mathcal C_{{\nu}}'(z)=\mathcal C_{{\nu}-1}(z)-\mathcal C_{{\nu}+1}(z)$ and $\mathcal C_{{\nu}-2}(z)+\mathcal C_{{\nu}}(z)=\frac{2({\nu}-1)}{z}\mathcal C_{{\nu}-1}(z)$, where $\mathcal C_{{\nu}}(z)$ stands both for $J_{\nu}(z)$ and $Y_{\nu}(z)$ (see \cite[\S 9]{abram}). This proves the lemma.
\endproof
\end{lem}

\begin{lem}\label{inducross}
The following identities hold
\begin{eqnarray}
Y_{\nu}(z)J_{\nu}^{(k)}(z)-J_{\nu}(z)Y_{\nu}^{(k)}(z)&=&\frac{2}{\pi z}\left(r_k+R_{{\nu},k}(z)\right),\label{num1}\\
Y_{\nu}'(z)J_{\nu}^{(k)}(z)-J_{\nu}'(z)Y_{\nu}^{(k)}(z)&=&\frac{2}{\pi z}\left(q_k+Q_{{\nu},k}(z)\right),\label{num2}
\end{eqnarray}
for all $k> 2$ and ${\nu}\geq 0$, where $r_k, q_k\in \{0,1, -1\}$, and $Q_{{\nu},k}(z)$, $R_{{\nu},k}(z)$ are finite sums of quotients of the form $\frac{c_{\nu,k}}{z^m}$, with $m\geq 1$ and $c_{\nu,k}$ a suitable constant, depending on ${\nu},k$.
\proof
We will prove (\ref{num1}) and (\ref{num2}) by induction. Identities (\ref{num1}) and (\ref{num2}) hold for $k=1$ and $k=2$ by Lemma \ref{case0}. Suppose now that
\small
\begin{eqnarray*}
Y_{\nu}(z)J_{\nu}^{(k)}(z)-J_{\nu}(z)Y_{\nu}^{(k)}(z)&=&\frac{2}{\pi z}\left(r_k+R_{{\nu},k}(z)\right),\\
Y_{\nu}'(z)J_{\nu}^{(k)}(z)-J_{\nu}'(z)Y^{(k)}(z)&=&\frac{2}{\pi z}\left(q_k+Q_{{\nu},k}(z)\right),
\end{eqnarray*}
\normalsize
hold for all ${\nu}\geq 0$. First consider 
\small
$$
Y_{\nu}'(z)J_{\nu}^{(k+1)}(z)-J_{\nu}'(z)Y_{\nu}^{(k+1)}(z).
$$
\normalsize
We use the recurrence relations $\mathcal C_{{\nu}+1}(z)+\mathcal C_{{\nu}-1}(z)=\frac{2{\nu}}{z}\mathcal C_{\nu}(z)$ and $2\mathcal C'(z)=\mathcal C_{\nu -1}(z)- \mathcal C_{\nu +1}(z)$, where $\mathcal C_{\nu}(z)$ stands both for $J_{\nu}(z)$ and $Y_{\nu}(z)$ (see \cite[\S 9]{abram}). We have
\footnotesize
\begin{multline}\label{laquzero}
Y_{\nu}'(z)J_{\nu}^{(k+1)}(z)-J_{\nu}'(z)Y_{\nu}^{(k+1)}(z)=Y_{\nu}'(z)(J_{\nu}')^{(k)}(z)-J_{\nu}'(z)(Y_{\nu}')^{(k)}(z)\\
=\frac{1}{4}\left[\left(Y_{{\nu}-1}(z)-Y_{{\nu}+1}(z)\right)\left(J_{{\nu}-1}(z)-J_{{\nu}+1}(z)\right)^{(k)}\right.\\
-\left.\left(J_{{\nu}-1}(z)-J_{{\nu}+1}(z)\right)\left(Y_{{\nu}-1}(z)-Y_{{\nu}+1}(z)\right)^{(k)}\right]\\
=\frac{1}{4}\left[\left(Y_{{\nu}-1}(z)J_{{\nu}-1}^{(k)}(z)-J_{{\nu}-1}(z)Y_{{\nu}-1}^{(k)}(z)\right)+\left(Y_{{\nu}+1}(z)J_{{\nu}+1}^{(k)}(z)-J_{{\nu}+1}(z)Y_{{\nu}+1}^{(k)}(z)\right)\right.\\
+\left.\left(J_{{\nu}+1}(z)Y_{{\nu}-1}^{(k)}(z)-Y_{{\nu}-1}(z)J_{{\nu}+1}^{(k)}(z)\right)+\left(J_{{\nu}-1}(z)Y_{{\nu}+1}^{(k)}(z)-Y_{{\nu}+1}(z)J_{{\nu}-1}^{(k)}(z)\right)\right]\\
=\frac{1}{4}\left[\frac{2}{\pi z}\left(r_k+R_{{\nu}-1,k}(z)+r_k+R_{{\nu}+1,k}(z)\right)\right.\\
+\left.\frac{2{\nu}}{z}\left(J_{\nu}(z)Y_{{\nu}-1}^{(k)}-Y_{\nu}(z)J_{{\nu}-1}^{(k)}(z)+J_{\nu}(z)Y_{{\nu}+1}^{(k)}(z)-Y_{\nu}(z)J_{{\nu}+1}^{(k)}(z)\right)\right.\\
-\left.\left(J_{{\nu}-1}(z)Y_{{\nu}-1}^{(k)}(z)-Y_{{\nu}-1}(z)J_{{\nu}-1}^{(k)}(z)+J_{{\nu}+1}(z)Y_{{\nu}+1}^{(k)}(z)-Y_{{\nu}+1}J_{{\nu}+1}^{(k)}(z)\right)\right]\\
=\frac{1}{4}\left[\frac{4}{\pi z}\left(2r_k+R_{{\nu}-1,k}(z)+R_{{\nu}+1,k}(z)\right)\right.\\
+\left.\frac{2{\nu}}{z}\left(J_{\nu}(z)\left(Y_{{\nu}-1}(z)+Y_{{\nu}+1}(z)\right)^{(k)}-Y_{\nu}(z)\left(J_{{\nu}-1}(z)+J_{{\nu}+1}(z)\right)^{(k)}\right)\right]\\
=\frac{1}{\pi z}\left(2r_k+R_{{\nu}-1,k}(z)+R_{{\nu}+1,k}(z)\right)\\
+\frac{{\nu}^2}{z}\left(J_{\nu}(z)\left(\frac{1}{z}Y_{\nu}(z)\right)^{(k)}-Y_{\nu}(z)\left(\frac{1}{z}J_{\nu}(z)\right)^{(k)}\right)\\
=\frac{2}{\pi z}\left[r_k+\frac{1}{2}\left(R_{{\nu}-1,k}(z)+R_{{\nu}+1,k}(z)\right)\right.\\
-\left.\frac{{\nu}^2}{z}\sum_{j=0}^k\frac{k!(-1)^{k-j}}{j!z^{k-j+1}}\left(r_j+R_{{\nu},j}(z)\right)\right].
\end{multline}
\normalsize

We prove now (\ref{num2})
\footnotesize
\begin{multline}\label{laqu}
Y_{\nu}(z)J_{\nu}^{(k+1)}(z)-J_{\nu}(z)Y_{\nu}^{(k+1)}(z)=\left(Y_{\nu}(z)J_{\nu}^{(k)}(z)-J_{\nu}(z)Y_{\nu}^{(k)}(z)\right)'\\
-\left(Y_{\nu}'(z)J_{\nu}^{(k)}(z)-J_{\nu}'(z)Y_{\nu}^{(k)}(z)\right)\\
=\frac{2}{\pi z}\left(-q_k-Q_{{\nu},k}(z)-\frac{r_k}{z}-\frac{R_{{\nu},k}(z)}{z}+R_{{\nu},k}'(z)\right).
\end{multline}
\normalsize
This concludes the proof.
\endproof
\end{lem}

\begin{corol}\label{formulaused}
The following formulas hold
\small
\begin{eqnarray*}
J_{\nu}(z)Y_{\nu}'''(z)-Y_{\nu}(z)J_{\nu}'''(z)&=&\frac{2}{\pi z}\left(\frac{2+{\nu}^2}{z^2}-1\right);\\
Y_{\nu}'(z)J_{\nu}'''(z)-J_{\nu}'(z)Y_{\nu}'''(z)&=&\frac{2}{\pi z^2}\left(1-\frac{3{\nu}^2}{z^2}\right);\\
Y_{\nu}'(z)J_{\nu}''''(z)-J_{\nu}'(z)Y_{\nu}''''(z)&=&\frac{2}{\pi z}\left(1-\frac{3+2{\nu}^2}{z^2}+\frac{{\nu}^4+11{\nu}^2}{z^4}\right).
\end{eqnarray*}
\normalsize
\proof
From Lemma \ref{inducross} (see in particular (\ref{laqu})) it follows
\small
\begin{multline*}
J_{\nu}(z)Y_{\nu}'''(z)-Y_{\nu}(z)J_{\nu}'''(z)=-\frac{2}{\pi z}\left[-q_2-Q_{{\nu},2}(z)-\frac{r_2}{z}-\frac{R_{{\nu},2}(z)}{z}+R_{{\nu},2}'(z)\right]\\
=\frac{2}{\pi z}\left(\frac{2+{\nu}^2}{z^2}-1\right).
\end{multline*}
\normalsize
Next we compute
\small
\begin{multline*}
Y_{\nu}'(z)J_{\nu}'''(z)-J_{\nu}'(z)Y_{\nu}'''(z)=\frac{2}{\pi z}\left[r_2+R_{{\nu},2}(z)-\frac{{\nu}^2}{z}\sum_{j=0}^2\frac{2(-1)^{2-j}}{j!z^{2-j+1}}\left(r_j+R_{{\nu},j}(z)\right)\right]\\
=\frac{2}{\pi z^2}\left(1-\frac{3{\nu}^2}{z^2}\right).
\end{multline*}
\normalsize
Finally, by (\ref{laquzero}) with $k=3$,  we have
\small
\begin{multline*}
Y_{\nu}'(z)J_{\nu}''''(z)-J_{\nu}'(z)Y_{\nu}''''(z)=\frac{2}{\pi z}\left[r_3+\frac{1}{2}\left(R_{{\nu}-1,3}(z)+R_{{\nu}+1,3}(z)\right)\right.\\
-\left.\frac{{\nu}^2}{z}\sum_{j=0}^3\frac{6(-1)^{3-j}}{j!z^{3-j+1}}\left(r_j+R_{{\nu},j}(z)\right)\right]\\
=\frac{2}{\pi z}\left(1-\frac{3+2{\nu}^2}{z^2}+\frac{{\nu}^4+11{\nu}^2}{z^4}\right).
\end{multline*}
\normalsize
\endproof
\end{corol}

{\bf Acknowledgments.} Large part of the computations in this paper have been performed by the second author in the frame of his PhD  Thesis under the guidance of 
the first author. 
The authors acknowledge financial support from the research project `Singular perturbation problems for differential operators',  Progetto di Ateneo of the University of Padova
and from the research project `INdAM GNAMPA Project 2015 - Un approccio funzionale analitico per problemi di perturbazione singolare e di omogeneizzazione'. 
The authors are members of the Gruppo Nazionale per l'Analisi Matematica, la Probabilit\`{a} e le loro Applicazioni (GNAMPA) of the Istituto Nazionale di Alta Matematica (INdAM).

\end{document}